\documentclass[12pt]{amsart}
\usepackage{amsthm,amsmath,amssymb}
\usepackage{addfont}
\usepackage{mathrsfs}
\usepackage[T1]{fontenc}
\newcommand{\n}{\noindent}
\numberwithin{equation}{section}

\def\ca{{\mathcal A}}
\def\cb{{\mathcal B}}

\def\bc{{\mathbb C}}

\def\bn{{\mathbb N}}

\def\a{\alpha}
\def\eps{\varepsilon}

\def\b{\beta}
\def\d{\delta}

\def\g{\gamma}
\def\k{\kappa}
\def\l{\lambda}       
\def\m{\mu}
\def\n{\nu}

\def\r{\rho}
\def\s{\sigma}
    
\def\t{\tau}

\def\f{\varphi}

\theoremstyle{plain}
\newtheorem{lemma}{Lemma}[section]
\newtheorem{proposition}[lemma]{Proposition}
\newtheorem{theorem}[lemma]{Theorem}
\newtheorem{corollary}[lemma]{Corollary}
\theoremstyle{definition}

\newtheorem{remark}[lemma]{Remark}
\newtheorem{definition}[lemma]{Definition}
\newtheorem{example}[lemma]{Example}

\begin{document}

\title[Jordan norms and  Grothendieck's inequalities]{\textsc{  Jordan norms for multilinear maps on  C*-algebras and Grothendieck's inequalities}}

\author[E.~Christensen]{Erik Christensen}
\address{\hskip-\parindent
Erik Christensen, Mathematics Institute, University of Copenhagen, Copenhagen, Denmark.}
\email{echris@math.ku.dk}
\date{\today}

\begin{abstract}   
There exists a generalization of the concept {\em completely bounded norm,} for multilinear maps on C*-algebras. We will use the word {\em Jordan norm,} for  this norm and denote it by  $\|.\|_J.$     The Jordan norm $\|\Phi\|_J$ of a multilinear map is obtained via  factorizations of  $\Phi$ in the form $$\Phi(a_1, \dots , a_n)  = T_0 \s_1(a_1)T_1 \dots T_{(n-1)}\s_n(a_n)T_n ,$$ where the maps $\s_i$ are  Jordan homomorphisms. We show that  any  bounded bilinear form on a pair of C*-algebras is Jordan bounded and satisfies $\|B\|_J \leq 2\|B\|. $  
 \end{abstract}

\maketitle

\section{Introduction}
Grothendieck's classical inequality may be presented in several ways and here we will focus on some aspects of the  extension of the following classical theorem theorem to the case of non commutative C*-algebras. The classical theorem we will focus on is the following, which may be found in  Pisier's survey article \cite{Pi3}, Theorem 2.1.

\begin{theorem} \label{ClasGr} 
There exists a positive real $K_G^\bc$ such that for any compact Hausdorff space $\Omega$ and any bounded bilinear form $B(f,g)$ on the continuous functions on $\Omega$ there exist probability measures $\m$ and $\n$ on $\Omega $ such that $$\forall f,g \in C(\Omega): \quad |B(f,g)|\leq K_G^\bc\|B\|\|f\|_{(2,\m)}\|g\|_{(2,\n)}.$$
\end{theorem} 
It is known that $1.33 < K_G^\bc < 1.41,$ \cite{Pi3}. 
 
This theorem was extended to the non commutative setting of C*-algebras  by Pisier in \cite{Pi4} under a mild approximation condition, which was removed by Haagerup in \cite{Ha2}, such that we now have the following non commutative version of Grothendieck's inequality, \cite{Pi4} equation (7.1).

\begin{theorem} \label{NcGr} 
Let $\ca$ be a C*-algebra and $B(x,y)$ a bounded bilinear form on $\ca, $ then there exist 4 states $\{\l, \k ,\m, \n\}$ on $\ca$ such that $$\forall x,y \in \ca:\quad |B(x,y)| \leq \|B\|\sqrt{\l(x^*x) + \n(xx^*)}\sqrt{\m(y^*y) + \n(yy^*)}.$$    
\end{theorem}

This theorem  proves the  validity of Theorem \ref{ClasGr} and  shows that   $K_G^\bc \leq2.$  

 We will not define and discus the concept named {\em completely bounded } here in  the introduction,  but give a short description of the facts we need in Section 2.   
In the article \cite{C3} we showed that the bilinear forms on a finite dimensional  abelian C*-algebra do have factorization properties which may be expressed naturally in the set up of operator spaces and completely bounded maps.  
It is well known that for a bounded bilinear form on a pair of  non commutative C*-algebras, one can not  expect that it is completely bounded in the way this concept is defined in \cite{CS} or \cite{Pa} and in Section 2.  On the other hand a bounded bilinear form on a C*-algebra $\ca$ may be considered to be a bounded operator from $\ca$ into its dual space $\ca^*$ and as such one may ask if it is completely bounded. 
This question was discussed and clarified by Pisier and Shlyakthenko in \cite{PS}, where they introduced the concept named {\em joint complete boundedness } for bilinear forms, which have the property that the associated operator from $\ca$ to the dual $\ca^*$ is completely bounded.  Not all bounded bilinear forms are joint completely bounded.  

Here we will try to get some insights into the structure of a general bounded bilinear form on a non commutative C*-algebra, and especially focus on an aspect which shows that also the non commutative Grothendieck inequality may be seen as a factorization result. To describe the non-commutative aspects we use in our approach, we will first make the remark that a bounded bilinear form on a unital  C*-algebra $\ca$ does not depend on the *-algebraic structure of $\ca,$ but only on the isometric structure of the Banach space $\ca.$
 Such a space is the non commutative analogue of the space $C(\Omega, \|.\|_\infty)$ for some compact Hausdorff space $\Omega$  and  the Banach-Stone Theorem says that in this case, the compact Hausdorff space  $\Omega$ is determined up to homeomorphism by the isometric  Banach space structure of $(C(\Omega), \|. \|_\infty) .$  In \cite{Ka} Kadison shows that for  non commutative C*-algebras the metric structure does not determine the *-algebraic structure, but any two C*-algebras which are isometrically isomorphic are isometrically isomorphic via a positive isometric mapping, which is a Jordan homomorphism. Further the book \cite{AS} by Alfsen and Schultz demonstrates that the compact convex set named the state space of a C*-algebra does not determine the C*-algebra, one also needs an orientation on the space.   This indicates that in the case of bounded bilinear forms you should look for an analogue of  the concept named {\em Stinespring representation,} but based on a Jordan representation instead of an algebraic *-representation. We show below, that this proposal works, when it is based  on the following definitions. 
 \begin{definition} \label{JrepDef} Let $\ca$ be a C*-algebra,  $G,H$  Hilbert spaces. 
 \begin{itemize} 
 \item[(i)] A self-adjoint map $\r: \ca \to B(H)$ is called a *-anti represen-tation if for all $x, y$ in $\ca,$ $\r(xy) = \r(y)\r(x).$  
 \item[(ii)] A positive self-adjoint  map $\s : \ca \to B(H)$ is said to be a Jordan representation if it is an orthogonal  sum of a *-representation and a *-anti representation. 
\item[(iii)]Let  $\ca_1, \dots, \ca_n $ be  C*-algebras.
An $n-$linear map $\Phi: \ca_1\times ... \times \ca_n  \to B(G,H) $ is said to be given by the 
Jordan-Stinespring representation \newline $\big( (K_1 , .., K_n) , ( \s_1 , .., \s_n) , ( T_0, .. , T_n)\big)$  if there exist Hilbert spaces $K_1, .. , K_n,$ operators $T_0 \in B(K_1, H), .., T_i \in B(K_{(i+1)}, K_i), . \newline ., T_n \in B(G, K_n)$ and Jordan representations $\s_i : \ca_i \to B(K_i)$ such that \begin{align*} \forall (a_1, \dots, a_n) \in  \ca_1 \times \dots \times \ca_n &:\\ \quad \Phi(a_1, \dots, a_n)& = T_0\s_1(a_1)T_1 \dots T_{(n-1)}\s_n(a_n) T_n.\end{align*}
\end{itemize}
\end{definition}  

If you wonder how a *-anti representation can occur, just think of the {\em prime example,} the transposition map $A \to A^t$ in $M_n(\bc).$ \medskip

Let $\Phi$ be a multilinear map on a product of C*-algebras with values in some $B(G, H).$ We  show in Section 3, that if $\Phi$ has a Jordan-Stinespring representation, then it is possible to define a norm, $\|\Phi\|_J,$  which we denote the Jordan norm by the definintion \begin{align} \label{Jnorm1}
&\|T\|_J := \inf\{ \|T_0\|\dots \|T_n\|\, : \\ \notag &\big((K_1, \dots, K_n), ( \s_1,\dots,  \s_n), (T_0, \dots, T_n)\big), \text{ Jor.-St. rep. of } \Phi \, \}. 
\end{align}

We can show that for any bounded bilinear form $B$   on a pair of  C*-algebras $\ca $ and $\cb,$ we have $\|B\| \leq \|B\|_J \leq 2 \|B\|,$ so the $J-$norm  behaves to a large extent as the completely bounded norm does for bilinear forms on commutative C*-algebras.   We also obtain a  similar result for the little Grothendieck inequality, which in the non commutative world of C*-algebras  may be formulated as follows. Let $T$ be a bounded operator from a C*-algebra $\ca$ into a Hilbert space $H,$ then  the Jordan norm $\|T\|_J$ satisfies $\|T \| \leq \|T\|_J \leq \sqrt{2} \|T\|.$  It is worth to remark that the reason, why we are able to get these results,  is primarily  Pisier's and Haagerup's non commutative Grotendieck inequalities. The last step, where it is put into the framework of Jordan representations, is inspired by Tomita's result from \cite{Ta}, or \cite{KR} Theorem 9.2.9, but here the actual use of that construction   comes nearly for free, as soon as it was detected as a missing link.   

 \section{On completely bounded maps} 
 The concept {\em complete boundedness}  has already been mentioned a couple of times, because it played a major role in our article  \cite{C3}, where we studied bilinear forms on  finite dimensional abelian C*-algebras. Since our scope here is to make the shift from the  abelian C*-algebras to the general non commutative C*-algebras as natural as possible, we will present the basic facts on {\em complete boundedness,}  and by an example show that it can not be used directly in the study of bounded bilinear forms on non commutative C*-algebras. On the other hand our replacement - the Jordan norm -   has inherited many properties, and it is obtained as a necessity in order to move the study of bilinear forms on $C(\Omega) $ to  the study of bilinear forms on general  C*-algebras. 
 
  Vern Paulsen has written the monograph \cite{Pa}, which presents the theory of completely bounded linear and multilinear maps on operator spaces  in  a most readable way. 
 
 \begin{definition} Let $\ca$ be a C*-algebra and $H$ a Hilbert space.
 \begin{itemize}
 \item[(i)] Let  $\f: \ca \to B(H)$ be a bounded linear map and let $n$ be a natural number. The map $\f_n : \ca \otimes M_n(\bc) \to B(H) \otimes M_n(\bc)$ is defined by $\f_n := \f\otimes\mathrm{id}_n.$ If $\sup\{\|\f_n\| \, : \, n \in \bn\} <\infty,$ then $\f$ is said to be  completely bounded and the supremum is defined as the completely bounded norm $\|\f\|_{cb} $ of $\f.$  
\item[(ii)]   
  Let $B$ be bounded bilinear form on a product of C*-algebras $\ca_1 \times  \ca_2.$ For any natural number $n$ the bilinear operator $$B_n : \big(\ca_1\otimes M_n(\bc)\big) \times \big(\ca_2 \otimes M_n(\bc)\big) \to M_n(\bc)$$ is defined by the equations \begin{align*}  \forall k, l \in \{1, \dots, n\} \,  &\forall  (X, \, Y) \,  \in M_n(\ca_1 ) \times  M_n(\ca_2):\\    B_n(X,Y)_{(k,l)} &= \sum_{j=1}^n B(X_{(k,j)}, Y_{(j,l)}).
  \end{align*}
   We say that $B$ is completely bounded if the set $\{\|B_n\|\, : \, n \in \bn \}$ is bounded, and if so the completely bounded norm $\|B\|_{cb}$ is defined as the supremum of this set, \cite{CS}. 
   \end{itemize}
  \end{definition} 
  
  Stinespring introduced in \cite{St} the concept named {\em completely positive. } A map $\f $ of  a C*-algebra into some other C*-algebra is defined to be completely positive if all the maps $\f_n:= \f\otimes \mathrm{id}_n$ are positive. He showed that a map $\f$ of a C*-algebra $\ca$ into $B(H)$ for some Hilbert space $H$ is completely positive if and only if it  can be decomposed in the following way. There exists a Hilbert space $K,$ a *-representation $\pi$ of $\ca$ on $K$  and a bounded operator $X $ in $B(H,K)$ such that $$\forall a \in \ca : \quad \f(a) = X^*\pi(a)X.$$  This theorem was later generalized to the setting of completely bounded linear or multilinear maps \cite{Pa} Theorem 8.4 and Corollary 17.7.  

The reason, why completely bounded linear or multilinear maps are nice to work with, is based on the fact, that they all have Stinespring representations, a concept we will define below. Moreover there always exist optimal Stinespring representations from where you can obtain the completely bounded norms directly.  
  
  \begin{definition} Let $\ca_1, \dots, \ca_n $ be  C*-algebras,  $H$ a Hilbert space and $\Phi: \ca_1 \times \dots \times \ca_n : \to B(H)$ a bounded $n$ linear map.  A Stinespring representation of $\Phi $ is a set $\big((K_1, \dots, K_n), ( \pi_1, \dots , \pi_n), (T_0,  \dots , T_n)\big) $ such that $K_i $ is a Hilbert space, $\pi_i $ is a *-representation of $\ca_i$ on $K_i$ and $T_i$ is  a bounded operator such that $T_0 $ is in $B(K_1, H),$ $T_i $ is in $B( K_{(i+1)}, K_i), $ $T_n$ in $B(H, K_n)$ such that  
  \begin{align*}  \forall (a_1, \dots, a_n) \in \ca_1 \times \dots \times \ca_n: & \\\quad \Phi(a_1, \dots, a_n)&  = T_0\pi_1(a_1)T_1 \dots T_{(n-1)}\pi_n(a_n)T_n.\end{align*}
  \end{definition} 
  It is elementary to see that  a bounded multilinear operator $\Phi$ with a Stinespring representation as above is completely bounded and the completely bounded norm  satisfies $\|\Phi\|_{cb} \leq \|T_0\| \dots \|T_n\|.$ It is  a theorem in this theory, that any completely bounded multilinear map $\Phi$ has a Stinespring representation as defined above such that $\|\Phi\|_{cb} = \|T_0\|\dots \|T_n\|.$

   We will give an example of  the differences between bilinear forms on abelian and bilinear forms on non commutative C*-algebras. First we recall that a bounded bilinear form $B$ on an abelian C*-algebra is automatically completely bounded and Grothendieck's constant $K_G^\bc$ is the smallest positive real which satisfies the inequality $\|B\|_{cb} \leq K_G^\bc\|B\|$ for all such bilinear forms, \cite{C3}. 
   
\begin{example} 
Let $(\d_n)_{( n \in \bn)}$  denote the canonical orthonormal  basis in $\ell^2(\bn).$ A bounded operator $X $ on $\ell^2(\bn)$ may then be represented by an infinite scalar matrix $(X_{(i,j)})$ where $i$ and $ j $ are natural numbers via the equation $X_{(i,j)} := \langle  X \d_j, \d_i \rangle.$ We can then define a bounded bilinear form $B$ on $B(\ell^2(\bn))$ by  $$B(X,Y) := (YX)_{(1,1)}. 
$$ It is obvious that $\|B\| = 1,$  but  $B$ is not completely bounded. To see this we will let $\{ e_{(i,j)} \, : \, i, j \in \bn\}$ denote the matrix units in $B(\ell^2(\bn))$ and we will fix a natural number $n. $ The matrix units in $ M_n(\bc) $  are denoted $f_{(i,j)}$ and we  can then define $X, Y $ in $B(\ell^2(\bn)) \otimes M_n(\bc)$ by the equations  $$ X = \sum_{k = 1}^n e_{(k,1)} \otimes f_{(1,k)} \text{ and } Y = \sum_{k = 1}^n e_{(1,k)} \otimes f_{(k,1)}.$$ It is not difficult to see that $X^* = Y$ and $X^*X = e_{(1,1) }\otimes I_n, $ so both $X$ and $Y$ are partial isometries with  $\|X\| = \|Y \| =1.$ On the other hand it follows from a computation that $$B_n(X,Y)_{(1,1)} = \sum_{k =1}^n B(X_{(1,k)},Y_{(k,1) })  = \sum_{k = 1}^n B(e_{(k,1)}, e_{(1,k)}) = n e_{(1,1)}, $$ so $\|B_n\| \geq n$ and $B$ is not completely bounded.
\end{example} 

\section{ The Jordan norms}
We introduced the function $\|.\|_J$ in  equation \eqref{Jnorm1}. It is not clear for which multilinear maps this function may be defined, and right away the definition does not seem to give a function which satisfies the triangle inequality.   In the abelian case it is clear that the completely bounded norm is a norm, and the fact that it is computable as a minimal value  over various factorizations is a theorem. Such an approach to the Jordan norm is not available here, but in this section we will  show that the equation \eqref{Jnorm1} actually defines a norm on multilinear operators which do have at least one Jordan-Stinespring representation.   
  In the case of a bilinear form $B$ on a product $\ca \times \cb$ of C*-algebras, or in the case of a bounded map $F$ of a C*-algebra into a Hilbert space,  we remind the reader, that in analogy with Dirac's notation we can, and will,  identify a vector $\xi$ in a Hilbert space $H$  with the operator $T_\xi :\bc \to H$ given by $T_\xi(z) := z \xi$ and similarly the adjoint  operator is given by $T_\xi^*(\eta) = \langle \eta, \xi \rangle.$ It is well known that these identifications between operators and vectors are isometric, but the map $\xi \to T^*_\xi$ is a conjugate linear isometry. 
  With these conventions in mind, the Jordan norms are also usable in the context of bilinear forms and for operators from a C*-algebra to a Hilbert space. 

\begin{definition} Let $\ca_1,\dots, \ca_n$ be C*-algebras, $G$ and  $H$  Hilbert spaces. We define the vector space $J(\ca_1, \dots, \ca_n,B(G,  H)) $ to be the space of all bounded multilinear maps of $\ca_1 \times \dots \times \ca_n$ into $B(G, H)$ which have a Jordan-Stinespring representation. \end{definition}

\begin{proposition}
Let $\ca_1, .., \ca_n$ be C*-algebras, $G$  and $H$  Hilbert \newline spaces. The positive function  $\|.\|_J $ on $J(\ca_1, \dots, \ca_n,B(G,  H)) $ defined by  \begin{align*} 
&\|\Phi\|_J := \inf\{\|T_0\|\dots \|T_n\| \, :\\ & \, \big((K_1, \dots, K_n), (\s_1, \dots,\s_n),( T_0, \dots, T_n)\big) \text{ Jor.-St. rep. of } \Phi\}, \end{align*} is a norm. 
\end{proposition} 

\begin{proof}
It is easy to see that $\|\Phi\| \leq \|\Phi\|_J,$ so  the only thing we need to prove is the triangle inequality for the function $\|.\|_J.$ 

Let $\Phi$ and $\Psi$ be in $J(\ca_1,  \dots , \ca_n, B(G, H)),$ and suppose that they  have the following Jordan-Stinespring representations,
\begin{align*} 
JS_\Phi& = \big((K_1, \dots, K_n), (\s_1, \dots,\s_n),( S_0, \dots, S_n)\big) \text{ for }  \Phi \\  JS_\Psi & = \big((L_1, \dots, L_n), \, (\t_1, \dots,\t_n),\, ( T_0, \dots, T_n)\big)
\text{ for } \Psi.  \end{align*} 
  We will also assume that the representations are normalized such that $\|S_0\| = \|S_n\|$  and $\|S_i\| = 1$ for $1 \leq i \leq n-1.$ Similarly we assume that $\|T_0 \| = \|T_n\|$ and $\|T_i\| = 1$ for $1 \leq i \leq n-1,$ and then \begin{equation} \label{normalize} 
\|S_0\| \dots \|S_n\| = \|S_0\|^2 \text{ and } \|T_0\| \dots \|T_n\| = \|T_0\|^2 
\end{equation} The following expressions defines a Jordan-Stinespring representation  $JS_{(\Phi+\Psi)}$ for $\Phi+ \Psi.$  \begin{align*} JS_{(\Phi+\Psi)}& :=  \big((K_1 \oplus L_1, \dots, K_n \oplus L_n),   (\s_1 \oplus \t_1 , \dots,\s_n \oplus \t_n), \\&  (\begin{pmatrix} S_0 & T_0\end{pmatrix}, \begin{pmatrix} S_1 & 0\\0 & T_1\end{pmatrix}, \dots,\begin{pmatrix} S_{(n-1)}  & 0\\0 & T_{(n-1)} \end{pmatrix},\begin{pmatrix} S_n\\T_n\end{pmatrix}) \big).
\end{align*} To see that the triangle inequality is valid, it is by \eqref{normalize} sufficient to prove the following inequality \begin{align*}   \|\begin{pmatrix} S_0 & T_0\end{pmatrix}\|&\| \begin{pmatrix} S_1 & 0\\0 & T_1\end{pmatrix}\| \dots \|\begin{pmatrix} S_{(n-1)}  & 0\\0 & T_{(n-1)} \end{pmatrix}\| \|\begin{pmatrix} S_n\\T_n\end{pmatrix}\|\\ &\leq \|T\|_0^2+\|S_0\|^2.
\end{align*}
This inequality follows from the fact that 
\begin{align*}
\|\begin{pmatrix}S_0 & T_0\end{pmatrix}\|^2 &= \| S_0S_0^* + T_0T_0^*\| \leq \|S_0\|^2 + \|T_0\|^2 ,\\
\|\begin{pmatrix} S_i  & 0\\0 & T_i \end{pmatrix}\| &= 1 \text{ for } 1 \leq i \leq n-1,\\
\|\begin{pmatrix} S_n\\T_n\end{pmatrix}\|^2 & = \| S_n^*S_n + T_n^*T_n\| \leq \|S_0\|^2 + \|T_0\|^2,
\end{align*} and the proposition follows. 
\end{proof}

It is interesting to know if the {\em infimum } used in the definition of the Jordan-Stinespring norm is  a {\em minimal value} ?  As mentioned above, this is so in the completely bounded case \cite{Pa}. We guess that the question has a positive answer. In the case  when the algebras $\ca_i$ are all finite dimensional and $H$ is finite dimensional too, we expect that a compactness argument based on norm topologies may be used to verify the guess in this finite dimensional setting.

\begin{remark} It seems clear that it may be possible to define other norms which most likely  are equivalent to the chosen Jordan norm. One suggestion is to replace the demand that the Jordan representations $\s_i$ used in the definition of a Jordan-Stinespring representation are asked to be a self-adjoint Jordan homomorphisms instead of an orthogonal sum of a *-representation and an *-anti representation. It is our guess, that this would give the same norm ?  \end{remark} 

\section{ The non commutative Grothendieck inequalities} 

The set up with Jordan bounded multilinear operators was primarily designed to give a way to express how the shift from the non commutative Grothendieck inequalities to the commutative ones may be expressed, such that the concepts used in the non commutative case become the concepts used in the commutative case, when the shift is made. We will show below that the non commutative Grothendieck inequality may be formulated in the following way: There exist two positive reals $K^J_G \leq 2 $ and $k_G^J \leq 2 $ such that any bounded bilinear  form $B$ on a pair of  C*-algebras  is Jordan bounded and satisfies  $\|B\|_J \leq  K^J_G\|B\|$ and  any bounded map $F$ of a C*-algebra into a Hilbert space is Jordan bounded and satisfies $\|F\|_J \leq \sqrt{k^J_G}\|F\|.$ 

In the commutative case we know that Grothendieck's inequalities may be formulated as follows \cite{C3}. The constants $k_G^\bc$ and $K_G^\bc$  are the smallest positive reals such that any bounded linear map $F$ of a commutative C*-algebra into a Hilbert space is completely bounded and satisfies $\|F\|_{cb} \leq \sqrt{k_G^\bc}\|F\|,$ and
any bounded bilinear form $B$  on a pair of  abelian C*-algebra is completely bounded and satisfies  $\|B\|_{cb} \leq K_G^\bc\|B\|. $

The proofs we give are  consequences of the works \cite{Pi4} and \cite{Ha2}, so we do not offer a new proof of the non commutative Grothendieck inequalities, but we suggest to look at the existing theorems in a slightly  different context, which - to us - seems to be a natural extension of our points of view on the classical inequalities as presented in \cite{C3} and \cite{C5}.  

The first  observation which makes it possible to use the Jordan norm in this context is the following proposition, which is inspired by the theory of Hilbert algebras, as used in the Tomita-Takesaki theory \cite{Ta} and \cite{KR} Theorem 9.2.9.  

\begin{proposition} \label{J}
Let $\ca$ be a C*-algebra and $\f$ a state on $\ca,$
then there exists a Hilbert space $H_\f,$ a unit vector $\xi_\f$ in $H_\f,$
a *-representation $\pi_\f$ of $\ca$ on $H_\f$ and a *-anti representation $\r_\f$ of $\ca$  on $H_\f$ such that $$ 
\forall a \in \ca: \f(a^*a) = \|\pi_\f(a)\xi_\f\|^2 \quad \text{ and } \quad \f(aa^*) = \|\r_\f(a)\xi_\f\|^2.$$
\end{proposition}

\begin{proof}
The Banach space double dual of $\ca$ may be represented on a Hilbert space as a von Neumann algebra,  and the state $\f$ extends naturally to a normal  state on this algebra. We let $\tilde{\ca} $ denote the sub C*-algebra of the double dual which is generated by $\ca$ and the unit in the double dual. We will use the GNS-construction \cite{KR} p. 279,   based on this state on $\tilde \ca $ to produce the set $\{ H_\f, \xi_\f, \pi_\f\}.$ To construct $\r_\f$ we supplement the unit vector $\xi_\f$ to an orthonormal basis $(\xi_\g)_{(\g \in \Gamma)}$ for $H_\f.$ We can then define a conjugate linear isometric involution $J$ on $H_\f$ by the equations $$ \forall \g \in \Gamma\, \forall z \in \bc: \quad J(z\xi_\g ) := \bar z \xi_\g. $$ We remind you, that in the case of a conjugate linear map, the adjoint operation on operators behaves differently, and in the case of the isometric involution $J$ we get the equation \begin{equation} \label{Jo} 
\forall \a, \b \in H_\f: \quad  \langle J\a, \b\rangle = \langle J\b , \a \rangle.
\end{equation} We can then define a *-anti representation $\r_\f$ by $$\r_\f(a) := J\pi_\f(a^*)J,$$ and it satisfies  
\begin{align*} \|\r_\f(a)\xi_\f\|^2  &=\langle \r_\f(a)^*\r_\f(a)\xi_\f, \xi_\f \rangle =  \langle \r_\f(aa^*)\xi_\f, \xi_\f \rangle\\ & = \langle J \pi_\f(aa^*)J\xi_\f, \xi_\f \rangle \text{ and by  } \eqref{Jo} \\ & = 
\langle J \xi_\f,  \pi_\f(aa^*)J\xi_\f \rangle \\  & = \langle  \xi_\f,  \pi_\f(aa^*)\xi_\f \rangle  =   \f(aa^*), 
\end{align*} and the proposition follows.

\end{proof}

\begin{theorem} \label{OJ} 
Let $\ca$ be a C*-algebra, $K$ a Hilbert space and $F: \ca \to K$ a bounded map. Then $F$ is Jordan bounded and satisfies  $\|F\|_J \leq \sqrt{2} \|F\|.$

\begin{proof}

By theorem 8.1 in \cite{Pi3} there exist states $\psi $ and $\f$ on $\ca$ such that $$\forall a \in \ca: \quad \|F(a)\|\leq  \|F\|\sqrt{ \psi(a^*a) + \f(aa^*)}.$$ Based on Proposition \ref{J} we construct Hilbert spaces $H_\psi$ with a unit vector $\xi_\psi,$ $H_\f$ with a unit vector $\xi_\f,$ a *-representation  $\pi_\psi$ of $\ca$ on $H_\psi$ and a *-anti representation $\r_\f$ on $H_\f$ such that $$ \forall a \in \ca: \quad \|\pi_\psi(a) \xi_\psi\|^2 = \psi(a^*a) \text{ and } \|\r_\f(a)\xi_\f\|^2 = \f(aa^*).$$
Then we see that there exists a bounded operator $S$ in $B(H_\psi \oplus H_\f, K)$ such that $\|S\|\leq \|F\| $ and $$ \forall a \in \ca : \quad F(a) = S \begin{pmatrix} \pi_\psi(a)& 0\\ 0& \r_\f(a)\end{pmatrix}\begin{pmatrix}\xi_\psi \\ \xi_\f\end{pmatrix} .$$ 

To get the connection of this to  a Jordan-Stinespring representation, we use the map $\xi \to T_\xi$ described in front of Definition 1.1, and we will look at the map $\ca  \ni a \to T_{F(a)}\in B(\bc,H).$   Let $\g$ denote the vector $( \xi_\psi , \xi_\f)$ in $H_\psi \oplus H_\f$ and let $T_\g$ in $B(\bc, H_\psi \oplus H_\f) $ be the associated operator, then $\|T_\g\| = \sqrt{2}.$  We have now obtained a Jordan-Stinespring representation $\big( (H_\psi \oplus H_\f), ( \pi_\psi \oplus \r_\f), (S, T_\g)\big)$ of the map $$\ca \ni a \to T_{F(a)} \in B(\bc, K),$$  which also shows that $\|F\|_J  \leq \sqrt{ 2}  \|F\|,$ and the theorem follows.
\end{proof}
\end{theorem}

The next theorem is proved  in the same way,  but this time from the well known non commutative Grothendieck inequality \cite{Pi4}, \cite{Ha2}.

\begin{theorem} \label{BJ}
Let $\ca, \,\cb$ be C*-algebras. Any bounded bilinear form  $B$  on $\ca \times \cb$ is Jordan bounded and $\|B\|_J \leq 2 \|B\|.$ 
\end{theorem} 

\begin{proof}
By Theorem 7.1 of \cite{Pi3} there exist 4 states $\{\l, \k, \m ,\n\}$ on $\ca$ such that \begin{equation} \label{GrIne}  \forall a \in \ca\, \forall b \in \cb: \,\, |B(a,b)| \leq \|B\| \sqrt{ \k(a^*a) + \l(aa^*)}\sqrt{ \m(b^*b) + \n(bb^*)}.\end{equation} 
We use the same constructions as above, so we produce  
4 Hilbert spaces with unit vectors $(H_\k, \xi_k), $ $(H_\l, \xi_\l),$ $(H_\m, \xi_\m),$
$(H_\n, \xi_\n),$  *-representions $\pi_\l$ of $\ca$  on $H_\l,$ $\pi_\m $ of $\cb$ on $H_\m$ and *-anti representations $\r_k$ of $\ca$ on $H_\k,$ $\r_\n$ of $\cb$ on $H_\n.$  The inequality \eqref{GrIne} shows that there exists a bounded operator $T$ in $B(H_\m \oplus H_\n , H_\l \oplus H_\k) $ such that $\|T\| \leq \|B\| $ and $T$ satisfies the following equation\begin{align} \label{FacB}
&\forall a \in \ca \, \forall b \in \cb: \\ \notag & B(a,b) = \langle T \begin{pmatrix} \pi_\m(b) & 0\\ 0 & \r_\n(b)\end{pmatrix}\begin{pmatrix} \xi_\m \\  \xi_\n\end{pmatrix} , \begin{pmatrix} \pi_\l(a^*) & 0
\\ 0 & \r_\k(a^*)\end{pmatrix}\begin{pmatrix} \xi_\l\\ \xi_\k\end{pmatrix} \rangle.
\end{align} 
 From here we can read a Jordan-Stinespring representation of $B$ and we get $$\|B \|_J \leq \|\begin{pmatrix} \xi_\l \\ \xi_\k\end{pmatrix}\| \|T\| \|\begin{pmatrix} \xi_\m \\ \xi_\n\end{pmatrix}\| = 2\|T\|,$$ and the theorem follows. 
\end{proof}

With this set up, it seems most likely  that the factorization we use in the proof of the theorem is not optimal with respect to the Jordan norm.  The reason being, that the factorizations we use all have  the property that $\|\xi_\m\| = \|\xi_\n\| = \|\xi_\k\| = \|\xi_\l\|.$ Right away we do not see any reason  why, an optimal factorization should satisfy these equations. 

The article \cite{PS} introduces the concept named {\em jointly completely bounded } for a bilinear form on a product of C*-algebras. This concept is based on the interpretation of a bounded  bilinear form as a bounded operator from one  C*-algebra into the dual of the other C*-algebra. Since the dual space of  a C*-algebra also is an operator space the concept of {\em complete boundedness} applies  to this operator too, and this is named {\em joint complete boundedness.} 
With the factorization given in \eqref{FacB} it is possible to express the joint completely bounded bilinear forms as those which may be factorized with a bounded operator $T$  in  diagonal form, which means, that with the notation from below, $T_{\l\n} =T_{\k\m} = 0.$  It seems to us that Haagerup and Musat's work on  {\em joint complete boundedness} in \cite{HM}, indicates that some difficult questions are still open.

It is immediate to obtain a splitting of a bilinear form $B = B_1+ \dots +B_4$ as presented in \cite{PS} at the bottom of page 186 in the following way. When $T$ is written as a block matrix 
$$T = \begin{pmatrix}T_{\l\m} & T_{\l\n}\\ T_{\k\m} & T_{\k\n}\end{pmatrix}, $$  then we may decompose it into a sum of 4 block matrices. This corresponds to the mentioned decomposition of $B,$ such that $$ B_1 \sim T_{\l\m} , \, B_2 \sim T_{\k\n}  , \, B_3 \sim 
T_{\l\n} , \,  B_4 \sim T_{\k\m} .$$

Based on the theorems \ref{OJ} and \ref{BJ} we make the following definition. 

\begin{definition} 
 The positive real $k^J_G$ is defined as the least real number such that for any bounded operator $F$ from a C*-algebra into a Hilbert space we have $\|F\|_J \leq \sqrt{k^J_G}\|F\|.$

The positive real $K_G^J$ is defined as the least real number such that for any bounded bilinear form $B$ on a pair of C*-algebras we have $\|B\|_J \leq K_G^J\|B\|.$   
\end{definition} 

It follows from the way the Jordan norms are defined that we always have $\|F\|_J \geq \|F\| $ and $\|B\|_J \geq \|B\|,$ so this and the  theorems \ref{OJ} and \ref{BJ} together imply that $$1 \leq k^J_G \leq 2 \text{ and } 1 \leq K^J_G \leq 2.$$ 

It will follow from  Theorem \ref{Pos} in the next section that $k^J_G \leq K_G^J.$

\section{Positive bilinear forms}

We find it interesting to take a special look at bounded positive bilinear forms on a C*-algebra, because, in the commutative case, the positive ones revealed some connections between the constants $k_G^\bc$ and $K_G^\bc,$ \cite{C5}. Part of this carries over to the non commutative case as we will show below. 

 We say that a bilinear form on a C*-algebra $\ca$ is positive, if for any $a$ in $\ca $ we have $B(a^*,a) \geq 0.$  
In the article Section 5 of \cite{C5} we  established the inequalities $\k_G^\bc \leq K_G^\bc \leq k_G^\bc/(2 - k_G^\bc) $ based on an investigation of positive bilinear forms. Unfortunately we are not able to follow the same path in the non commutative case, but we still get a link between the constants $k_G^J $ and $K_G^J $ which is based on a closer look at bounded positive bilinear forms on general C*-algebras.

A positive bilinear form $B$ on a C*-algebra $\ca$ induces a positive sesqui-linear form $[., .]_B $ on $\ca$ by the definition $$[x,y]_B:= B(y^*,x).$$ 
The kernel, defined as  $\{ x\in \ca : B(x^*,x) = 0\},$ may be non trivial, in which case we will work with cosets instead of operators and then we get a Hilbert space $K_B$  and a bounded map $ F_B : \ca \to K_B$ such that $$\langle F_B(x), F_B( y ) \rangle = [x, y]_B = B(y^*,x).$$ 
This means that $F_B$ is bounded and $\|F_B\|^2  \leq \|B\|.$ On the other hand the Cauchy-Schwarz inequality shows that $\|B\| \leq \|F_B\|^2$  so $\|F_B\|^2  = \|B\|.$  
If for some Hilbert space $L$ and some bounded operator $G:\ca \to L$ we have $\langle G(x), G(y) \rangle = B(y^*,x),$ 
then it is well known, that there exists an isometry $W$ of $K_B$ onto the range space of $G$ in $L$ such that  $G = WF_B.$ 
We have now set the stage to prove the following theorem. 

\begin{theorem} \label{Pos} Let $B$ be a bounded positive bilinear form on a C*-algebra $\ca$ and $F_B:\ca  \to K_B$ the bounded map of $\ca$ into  the Hilbert space associated to $B.$  Then  $\|B\| = \|F_B\|^2,$  $\|B\|_J = \|F_B\|_J^2$  and the constant $k^J_G$ is the smallest positive real $c$ such that for any C*-algebra $\ca$ and any bounded positive bilinear form $B$ on $\ca,$ $\|B\|_J \leq c \|B\|.$ \end{theorem}
\begin{proof} We showed above that $\|B\| = \|F_B\|^2.$ 

 To  $\eps >0$ there exists  a Jordan-Stinespring representation of $F_B$  \begin{equation}  \label{JF}  \forall a \in \ca : \quad F(a) = T \s(a) T_\xi, \, \|T\|\|\xi\| \leq \|F_B\|_J + \eps .\end{equation}  
We can then obtain a Jordan-Stinespring representation for $B$ by 
\begin{align}  \label{JB1} &\forall a,b  \in \ca: \\& \notag   B(a,b)=  T_\xi^*\s(a) (T^* T) \s(b)  T_\xi , \, \|\xi\|\|T^*T\|\|\xi\| \leq (\|F_B\|_J + \eps)^2 .\end{align}  
Then we get $\|B\|_J \leq \|F_B\|_J^2.$ 

On the other hand let us consider a Jordan-Stinespring representation of $B$ such that \begin{align}  \label{JB2} &\forall a,b  \in \ca: \\ \notag & B(a,b)=  T_\eta^* \s_l(a) T \s_r(b)  T_\xi , \, \|\eta\|\|T\|\|\xi\| \leq \|B\|_J + \eps .\end{align}
From here we can construct a self-adjoint positive Jordan-Stinespring representation.  To do this, we introduce some notation. First we assume that $\eta $ and $\xi$ are unit vectors in some Hilbert spaces $L$ and $R.$ Then $\s_l$ is an orthogonal sum of a *-representation and a *-anti representation of $\ca $ on $L$ and $\s_r$ is an orthogonal sum of a *-representation and a *-anti representation of $\ca $ on $R$. Finally $T$ is a bounded operator from $R$ to $L$ such that $\|T\| \leq \|B\|_J + \eps.$ We define  $\g = 2^{-(1/2)} (\eta, \xi)$ as  a unit vector in $L \oplus R$ and  $\s:= \s_l \oplus \s_r$ as an orthogonal sum of a *-representation and a *-anti representation of $\ca$ on $L \oplus R. $  Finally we define a bounded self-adjoint operator $S$ on $L \oplus R$ such that $\|S\| = \|T\|$  via the block matrix $$S := \begin{pmatrix}
o& T \\ T^* & 0
\end{pmatrix},$$ and we have obtained a self-adjoint Jordan-Stinespring representation of $B$ in the following way 

\begin{align*}
\label{JB3} \forall a,b  \in \ca: &\\   \quad  T_\g^* \s(a) S \s(b)  T_\g &= \frac{1}{2} \langle \s_l(a)T\s_r(b)\xi, \eta\rangle + \frac{1}{2} \langle \s_r(a)T^*\s_l(b) \eta , \xi \rangle  \\
&= \frac{1}{2} \langle \s_l(a)T\s_r(b)\xi, \eta\rangle + \frac{1}{2} \langle \eta, \s_l(b^*)T \s_r(a^*) \xi  \rangle  \\
&= \frac{1}{2} B(a,b)  +  \frac{1}{2} \overline{B(b^*,a^*)}  \\ &= \frac{1}{2} [b,a^*]_B   +  \frac{1}{2} \overline{[a^*,b]_B}  \\& = [b,a^*]_B = B(a,b), \text{ and } \\  \|\g\|\|S\|\|\g\| &\leq \|B\|_J + \eps .\end{align*}

Let us define $Q$ as the orthogonal projection of $L \oplus R$ onto the closure of the linear subspace of $L \oplus R$ 
given as $\{\s(b)\g\, : \, b \in \ca\},$ then the positivity of $B$ implies that $QSQ $ is a positive bounded operator with $\|QSQ\| \leq \|B\|_J + \eps,$ and we have obtained a positive Jordan-Stinespring representation of $B$ and in turn obtained a bounded operator $G: \ca \to Q(L\oplus R)$ with a Jordan-Stinespring representation given by 
\begin{equation} \label{GjB}  \forall a \in \ca: \quad G(a):= (QSQ)^{(1/2)} \s(a) T_\g \text{ and } \|G\|_J \leq \sqrt{\|B\|_J + \eps}.\end{equation} Then $$\langle G(b), G(a^*)\rangle =  B(a,b)= \langle F_B(b), F_B(a^*)\rangle,$$  and there exists an isometry $W$ of $K_B$ onto $G(L \oplus R) $ such that for any $a$ in $\ca$ we have $G(a) = WF_B(a)$ and in particular we have a Jordan-Stinespring representation of $F_B$ given by 
 \begin{equation} \label{GjFB}   \quad F_B(a):= W^*(QSQ)^{(1/2)} \s(a) T_\g,\,  \,   \|W^*(QSQ)^{(1/2)} \|\|\g\| \leq \sqrt{\|B\|_J+ \eps},\end{equation} and $\|F_B\|_J \leq \sqrt{\|B\|_J} .$ The theorem follows. \end{proof}

The theorem has the following immediate corollary.

\begin{corollary} The    Jordan-Grothendieck constants satisfy $$1.27 < k_G^\bc \leq  k^J_G \leq K_G^J \leq 2 \text{ and } 1.33 < K_G^\bc \leq K_G^J.$$ 
\end{corollary}  
\begin{proof} For a bounded linear map of an abelian  C*-algebra into a Hilbert space, the completely bounded norm and the Jordan norm agree,  so $k_G^\bc \leq k_G^J.$ 
By \cite{Pi3} Theorem 5.1 we know that $k_G^\bc = \frac{4}{\pi} > 1.27 , $ so the first string of  inequalities  follows from this and the theorem. From \cite{Pi3},  a little after equation (4.2), it follows that $K_G^\bc > 1.33. $ Since the completely bounded norm and the Jordan norm agree in the commutative case we get $1.33 < K_G^\bc \leq K_G^J.$  
\end{proof}

\end{document}